\documentclass{amsart}

\usepackage{graphicx,amsmath,amsthm,amsopn,amstext,amscd,amsfonts,amssymb}

\def\Today{\number\day\ \ifcase\month\or January\or
February\or  March\or April\or May\or June\or July\or
August\or  September\or October\or November\or December\fi\
\number\year}

\newcount\mins \newcount\hours \hours=\time \mins=\time
\divide\hours by60 \multiply\hours by 60 \advance\mins
by-\hours
\divide\hours by60

\def\now{ \ifnum\hours>11 \ifnum\hours>12 \advance\hours by
-12 \fi
\number\hours:\ifnum\mins<10 0\fi \number\mins\ pm,\ \else
\ifnum\hours=0 \hours=12 \fi
\number\hours:\ifnum\mins<10 0\fi \number\mins\ am,\ \fi}

\usepackage{fancybox}

\def\fig#1{
\includegraphics[height=2.5in]{#1.eps}
}

\newcommand{\wo}{\mathbb{E}}                                                    
\newcommand{\wwo}[2]{\mathbb{E}\left(\left.{#1}\right|{#2}\right)}              
\DeclareMathOperator{\Var}{Var}                                                 
\newcommand{\wVar}[2]{\Var\left(\left.{#1}\right|{#2}\right)}                   

\def\<{\langle}
\def\>{\rangle}

\def\<{\langle}
\def\>{\rangle}

\newcommand{\rf}[1]{(\ref{#1})}

\newcommand{\calF}{{\mathcal F}}

\newcommand{\calL}{{\mathcal L}}

\newcommand{\eps}{\varepsilon}
\newcommand{\la}{\lambda}

\newcommand{\be}{\begin{equation}}
\newcommand{\ee}{\end{equation}}

      \newtheorem{theorem}{Theorem}[section]
       \newtheorem{proposition}[theorem]{Proposition}

\theoremstyle{remark}
       
\theoremstyle{definition}

\def\<{\langle}
\def\>{\rangle}

\def\<{\langle}
\def\>{\rangle}

\newcommand{\Fsu}{{\mathcal{F}}_{s,u}}

\title{Classical bi-Poisson process: an invertible quadratic harness}

\author{
W{\l}odzimierz  Bryc
}
\thanks{\noindent Research partially supported by NSF
grants \#INT-0332062, \#DMS-0504198, and by the C.P. Taft Memorial Fund.}
\address{
Department of Mathematics,
University of Cincinnati,
PO Box 210025,
Cincinnati, OH 45221--0025, USA}
\email{Wlodzimierz.Bryc@UC.edu}

\author{Jacek Weso{\l}owski}
\address{ Faculty of Mathematics and Information Science\\
Warsaw University of Technology\\ pl. Politechniki 1\\ 00-661
Warszawa, Poland}
\email{wesolo@alpha.mini.pw.edu.pl}

\date{August 18,  2005.  \\ {\tt Printed \now\today\  File: \jobname.TEX}}


\begin{document}

\maketitle
\begin{abstract} We give an elementary construction of
a time-invertible Markov process which is discrete except at one instance.
The process is one of the quadratic harnesses studied in
\cite{Bryc-Wesolowski-03}, \cite{Bryc-Matysiak-Wesolowski-04},  and
\cite{Bryc-Matysiak-Wesolowski-04b}. It can be regarded as a random
joint of two independent Poisson processes.
\end{abstract}

\section{Introduction}
According to \cite{Watanabe-74}, a stochastic process $(X_t)_{t>0}$ has the time inversion property, if it has the
same finite-dimensional distributions as the process $(t
X_{1/t})_{t>0}$. In 
papers \cite{Gallardo-Yor04} and \cite{Lawi-05} the authors
give criteria for the time-invertibility of Markov processes with transition probabilities that have smooth densities with respect to the
Lebesgue measure.

In this note we give a new elementary example
of a time-invertible Markov process for which  all transitions except to time $t=1$ are discrete,
see Proposition \ref{Time Inversion}.
This improves upon \cite[Corollary 3.4]{Bryc-Wesolowski-04}, where we gave a less elementary example of a time-invertible Markov
 process that had transition probabilities with a discrete component.
Both examples are  particular cases of a more general family of Markov processes which  in \cite{Bryc-Matysiak-Wesolowski-04} we
called the bi-Poisson processes.
According to \cite[Example 4.8 and Proposition 4.13]{Bryc-Matysiak-Wesolowski-04},
  a bi-Poisson process with parameters $(\eta,\theta, q)$  is a square-integrable Markov
 process $(X_t)$ which is uniquely determined by the following three properties:
\begin{equation}\label{cov}
\wo(X_t)=0,\: \wo(X_tX_s)=\min\{t,s\},
\end{equation}
\begin{equation}
\label{LH}
\wwo{X_t}{\Fsu}=\frac{u-t}{u-s} X_s+\frac{t-s}{u-s} X_u,
\end{equation}
\begin{multline} \label{q-Var}
\wVar{X_t}{\Fsu} =   \\
  \frac{(u-t)(t-s)}{u-qs} \left( 1+ \eta\frac{uX_s-sX_u}{u-s} +\theta \frac{X_u-X_s}{u-s}
-(1-q)\frac{(uX_s-sX_u)(X_u-X_s)}{(u-s)^2} \right),
\end{multline}
for all $0\leq s<t<u$,
where
$$
\Fsu=\sigma\{X_t: 0\leq t \leq s  \mbox{ or } t>u\}.
$$
Property \rf{LH} says that the bi-Poisson process is a harness,  see
\cite{Mansuy-Yor04}. Condition \rf{q-Var} means that it is a
quadratic harnesses, see \cite{Bryc-Matysiak-Wesolowski-04}. The
adjective "classical"   refers to the value of parameter $q=1$,
compare \cite[Section 4.2]{Bryc-Matysiak-Wesolowski-04} and
\cite[Section 4.2]{Bryc-Wesolowski-03}. In Proposition \ref{two
Poissons} we show that $(X_t)$ can be constructed by joining
together  two independent Poisson processes with the same random
gamma intensity.  This is accomplished by appropriate affine
transformations and deterministic changes of time.

In \cite{Bryc-Wesolowski-04} we use orthogonal polynomials to construct the
 transition probabilities of the bi-Poisson process when $q=0$, and we show that its
 univariate distributions form a semigroup with
 respect to a certain generalized convolution related to free probability.
 The univariate distributions of the general bi-Poisson process
 were implicitly identified in \cite[Example 4.8]{Bryc-Matysiak-Wesolowski-04}, and the corresponding Markov process
 is under construction in
\cite{Bryc-Matysiak-Wesolowski-04b}.
This  construction relies heavily on  cumbersome identities between certain multi-parameter families
 of orthogonal polynomials, and  identifies the transition probabilities in implicit form only.
However, when $q=1$, the explicit transitions probabilities can be read out.
 They turn out to be
 related to the pure birth and the pure death processes and are amenable  to
  explicit elementary  analysis.

In this paper we present an elementary construction of the
bi-Poisson process with parameters $(\eta,\theta,1)$. Throughout
most of the paper the value of the third parameter is fixed as $q=1$, in which case we  say that we consider a bi-Poisson process with
parameters $(\eta,\theta)$, skipping
the third parameter of the triple.


The plan of the paper is as follows. In Section \ref{S: Constr} we give the construction of the process.
In Section \ref{S: Moments} we verify that the construction indeed gives a bi-Poisson process.
In Section \ref{S: AddProp} we deduce some additional properties, including
time-invertibility.

\section{Construction}\label{S: Constr}
 It is known, see \cite[Proposition 4.13]{Bryc-Matysiak-Wesolowski-04} that a bi-Poisson process $(X_t)$
with parameters $(\eta,\theta)$ satisfies $\eta\theta\geq0$.

In the degenerate case $\theta\eta=0$ it is known that the bi-Poisson process $(X_t)$ is either
$X_t=B_t$, where $(B_t)$ is the standard Brownian motion for
$\eta=\theta=0$,  or $X_t=\theta N_{t/\theta^2}-t/\theta$, where
$(N_t)$ is the standard Poisson process when  $\eta=0$, $\theta\ne 0$, see \cite[Theorem 1]{Wesolowski93}.
Passing to the time inverse $(tX_{1/t})$, we  see
that
$X_t=\eta t N_{1/(t\eta^2)}-1/\eta$ in the remaining degenerate case $\theta=0$, $\eta\ne 0$.

We will  therefore concentrate on the case $\theta\eta> 0$. Passing
to $(-X_t)$ preserves \rf{q-Var} replacing parameters
$(\eta,\theta)$ by $(-\eta,-\theta)$, so we may assume
$\eta,\theta>0$. Replacing  process $(X_t)$ by  process
$(\sqrt{\eta/\theta}X_{t\theta/\eta})$, we get the bi-Poisson
process with parameters $(\sqrt{\eta\theta},\sqrt{\eta\theta})$.
Thus without loss of generality we may  assume that $\eta=\theta>0$.

\begin{figure}[htb]
\fig{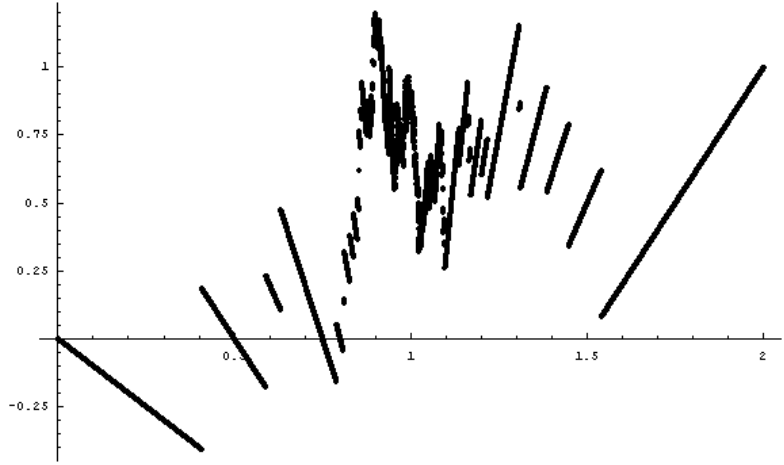} \caption{Simulated sample trajectory of the bi-Poisson
process. \label{Fig3} The process follows the segments $\ell_j: y=
\theta (1-t)j-t/\theta$, $0<t<1$, $j=0,1,2,\dots$ except for the
upwards jumps, and then follows the half-lines $y=\theta (t-1)j -
1/\theta$, $t>1$, $j=\dots, 1,0$  with the downwards jumps.}
\end{figure}

The moment of time $t=1$ is preserved by the time-inversion and plays a special role in the construction.
The bi-Poisson process traverses  a family of deterministic lines, with jumps in the upwards direction when $t<1$
and in the downwards direction when $t>1$, see Fig. \ref{Fig3}.

The process is determined by specifying
an integer that describes the line being followed at time $t$. The
integers that describe the upwards jumps form a  linear
pure birth process with immigration with the time   transformed to run on the interval $[0,1)$.
 At time $t=1$ instead of being infinite, the process takes
the continuous  spectrum of real values. For $t>1$, the downwards
jumps form a linear pure death process which "returns from $\infty$"
by a Poisson entrance law, again with the time transformed to run on
the interval $(1,\infty)$. The deterministic time transformations
are logarithmic and introduce a rather simple  non-homogeneity into
the birth rates and the death rates of the process. However, they
force infinite number of jumps before and after $t=1$.

For a more formal description of $(X_t)$, we set
\begin{equation}\label{X2Z}
X_t=\begin{cases}
\theta(1-t)Z_{t} - \frac{t}{\theta},& 0\leq t<1\;,
\\
 \theta Z_1 -\frac{1}{\theta}, & t=1,
\\
\theta( t-1)Z_t -\frac{1}{\theta},& t>1,
\end{cases}
\end{equation}
where random variables $Z_t$ are $\{0,1,2,\dots\}$-valued for $t\ne 1$.
We will construct the appropriate process $(Z_t)_{t\geq 0}$ in three steps: we first define $(Z_t)_{0\leq t<1}$ as
 a pure birth process, then we extend it to $t=1$ by passing to the limit, and finally we extend the process
 to $t>1$ as a pure death process with a $Z_1$-dependent Poisson entrance law.

\subsection{The pure birth phase}\label{The pure birth phase}
As $(Z_t)_{0\leq t<1}$ we take the non-homogeneous linear pure birth
process with immigration with the birth rate
$$\la_n(t)=\frac{n+\frac{1}{\theta^2}}{1-t}.$$
 The properties of such a process are well known. For small enough $|z|$,
\cite[Exercise 5.1]{Parzen-62} gives the following generating function of the
transition probabilities for the more general non-homogeneous linear pure birth process with the birth rate
 $\la_n(t)=\nu(t)+n\la(t)$.
$$
\sum_{k=0}^\infty z^k p_{j, j+k}(s,t)=z^{\nu(s)/\la(s)-\nu(t)/\la(t)}
 \left(\frac{ p(s,t)}{1-z(1-p(s,t))}\right)^{j+\nu(s)/\la(s)},
$$
where $p(s,t)=e^{-\int_s^t \la(u)du}$. In our setting,  $p(s,t)=\frac{1-t}{1-s}$ and with $z=e^{u(1-t)}$ we get
\begin{equation}\label{MGF}
\wwo{e^{u(1-t)(Z_t-Z_s)}}{\calF_{\leq s}}=\left(\frac{1-t}{1-s-(t-s)e^{u(1-t)}}\right)^{Z_s+1/\theta^2},
\end{equation}
where $\calF_{\leq s}=\sigma\{X_r: r\leq s\}=\sigma\{Z_r: r\leq
s\}$. Thus the conditional distribution $\calL(Z_t-Z_s|Z_s)$ is
negative binomial with parameters $r=Z_s+1/\theta^2$,
$p=(1-t)/(1-s)$.
 (Table \ref{Distributions}  lists the parameterizations of the distributions we use in this note.)

\begin{table}[hbt]
\begin{tabular}{||c|c|c|c||} \hline\hline
Name  & Parameters & Distribution & $E(e^{uZ})$ \\ \hline\hline
Poisson& $\la>0$ & $e^{-\la}\la^k/k!,k=0,1,2,\dots$&
$\exp(\la(e^u-1))$\\
Gamma & $p>0$, $\sigma>0$&
$f(z)=\frac{1}{\sigma^p\Gamma(p)}x^{p-1}e^{-x/\sigma}$
&$(1-\sigma u)^{-p}$  \\
Negative Binomial&$r>0$, $0<p<1$& $\frac{\Gamma(k+r)}{\Gamma(r)k!}p^r (1-p)^{k}, k=0,1,\dots$
&$\frac{p^r}{(1-(1-p)e^u)^r}$ \\
Binomial & $n\geq 0$, $0\leq p \leq 1$ & $(^n_k)p^k(1-p)^{n-k}$,
 $k=0,1,\dots,n$.
&$(1-p+pe^u)^n$\\
\hline \hline
\end{tabular}
\caption{Laws that appear as  transition probabilities.\label{Distributions}}
\end{table}

Differentiating \rf{MGF} at $u=0$ and using \rf{X2Z} we verify that  $(X_t,\calF_{\leq t})_{0\leq t<1}$ is a martingale.
Setting $s=0$, from \rf{MGF} we get $\wo(X_t)=0$ and differentiating \rf{MGF} again, after a calculation we get $\wo(X_t^2)=t$.

In particular, $X_1=\lim_{t\to 1}X_t$ converges almost surely, as
does $Z_1=\lim_{t\to 1}(1-t)Z_t$. Taking the limit in \rf{MGF} we
see that
\begin{equation}\label{MG Gamma}
\wwo{e^{uZ_1}}{Z_s}=\left(1-u(1-s)\right)^{-Z_s-1/\theta^2},
\end{equation} thus $\calL(Z_1|Z_s)$ is gamma with shape parameter
$p=Z_s+\frac{1}{\theta^2}$ and scale parameter $\sigma=1-s$, see
Table \ref{Distributions}. In particular, $Z_1$ is gamma with
$r=1/\theta^2$, $\sigma=1$,  and the support of $X_1$ is
$[-1/\theta,\infty)$.

\subsection{The pure death phase}
We now extend $(Z_t)$ from $0\leq t\leq 1$ to $t>1$ by specifying
$(Z_t)_{t>1}$ as a pure death process with the death rate
$$\mu_n(t)=\frac{n}{t-1},$$
and the $Z_1$-dependent entrance law $\calL(Z_t|Z_1)$ which we take as the Poisson
law with parameter $\la={Z_1}/{(t-1)}$. Thus
\begin{equation}\label{MG Poiss}
\wwo{e^{uZ_t}}{Z_1}=\exp\left(Z_1\frac{e^u-1}{t-1}\right).
\end{equation}

A well known property of the linear pure birth process is that  for $1<s<t$ the transition probabilities
$\calL(Z_t|Z_s)$ are binomial with parameters $n=Z_s$,
$p=(s-1)/(t-1)$, so
\begin{equation}\label{MG Bin}
\wwo{e^{uZ_t}}{Z_s}=\left(\frac{t-s+(s-1)e^u}{t-1}\right)^{Z_s}.
\end{equation}

 The Poisson distribution is indeed the entrance law: given $1<s<t$
we have
$$\Pr(Z_t=i|Z_1)=\sum_{n=i}^\infty \Pr(Z_t=i|Z_s=n)\Pr(Z_s=n|Z_1).$$
Indeed, the right hand side is
$$\frac{Z_1^i}{ i!  (t-1)^i} e^{-Z_1/(s-1)}
\sum_{n=i}^\infty \frac{Z_1^{n-i}(t-s)^{n-i}}{(n-i)!(s-1)^{n-i}(t-1)^{n-i} } =
\frac{Z_1^i}{i!  (t-1)^i}e^{-{Z_1}/{(t-1)}}.
$$

We now verify that the two pieces of the process fit together into  a well defined Markov process $(Z_t)_{t\geq 0}$.
For $0\leq s<1<t$, by conditioning on $Z_1$ we get
$$\Pr(Z_t=j|Z_s)=
\frac{1}{j!(t-1)^j}\wwo{Z_1^je^{-Z_1/(t-1)}}{Z_s}$$ $$ =
\int_0^\infty \frac{ x^{Z_s+j+1/\theta^2-1}
}{j!(t-1)^{j}(1-s)^{Z_s+1/\theta^2}
\Gamma(Z_s+1/\theta^2)}e^{-{x}/{(1-s)}} dx
$$
$$
=\frac{\Gamma(i+j+1/\theta^2)}{j!\Gamma(i+1/\theta^2)}\left(\frac{t-1}{t-s}\right)^{i+1/\theta^2}\left(\frac{1-s}{t-s}\right)^j.
$$
Thus $\calL(Z_t|Z_s)$ is negative binomial with $r=Z_s+1/\theta^2$ and $p=(t-1)/(t-s)$. In particular,
 $Z_t$ is negative binomial with $r=1/\theta^2$ and $p=1-1/t$. An elementary calculation shows that
 $X_t$ defined by \rf{X2Z} has mean zero and variance $t$.

 A straightforward calculation leads now to the verification of the Chapman-Kolmogorov equations in the remaining two cases:
\begin{enumerate}
\item If $0<s_1<s_2<1<t$ and $i,k\geq 0$ then
$$
\Pr(Z_t=k|Z_{s_1}=i)=\sum_{j=i}^\infty \Pr(Z_t=k|Z_{s_2}=j)\Pr(Z_{s_2}=j|Z_{s_1}=i).
$$
Indeed, the right hand side is
$$
\frac{1}{k!\Gamma(i+1/\theta^2)}\left(\frac{t-1}{1-s_1}\right)^{i+1/\theta^2}
\sum_{j=i}^\infty\frac{\Gamma(j+k+1/\theta^2)(1-s_2)^{i+k+1/\theta^2}(t-1)^{j-i}(s_2-s_1)^{j-i}}{(j-i)! (1-s_1)^{j-i}
(t-s_2)^{j-i}(t-s_2)^{i+k+1/\theta^2}}
$$
$$
=\frac{\Gamma(i+k+1/\theta^2)}{k!\Gamma(i+1/\theta^2)}\left(\frac{t-1}{1-s_1}\right)^{i+1/\theta^2}
\frac{(1-s_2)^{i+k+1/\theta^2}}{(t-s_2)^{i+k+1/\theta^2}} \left(1-\frac{(t-1)(s_2-s_1)}{(1-s_1)(t-s_2)}\right)^{-(i+k+1/\theta^2)}
$$
$$
=\frac{\Gamma(i+k+1/\theta^2)}{k!\Gamma(i+1/\theta^2)}\left(\frac{t-1}{t-s_1}\right)^{i+1/\theta^2}\left(\frac{1-s_1}{t-s_1}\right)^k
.$$

\item If $0<s <1<t_1<t_2$ and $i,k\geq 0$ then
$$
\Pr(Z_{t_2}=k|Z_{s}=i)=\sum_{n=k}^\infty \Pr(Z_{t_t}=k|Z_{t_1}=n)\Pr(Z_{t_1}=n|Z_{s}=i).
$$
Indeed, the right hand side is
$$
\frac{(1-s)^k(t_1-1)^{i+k+1/\theta^2}}{k!\Gamma(i+1/\theta^2)(t_1-s)^{i+k+1/\theta^2}(t_2-1)^k }\sum_{n=k}^\infty\frac{\Gamma(n+i+1/\theta^2)}{(n-k)!}
\left(\frac{(t_2-t_1)(1-s)}{(t_2-1)(t_1-s)}\right)^{n-k}
$$
$$
=\frac{\Gamma(k+i+1/\theta^2)}{k!\Gamma(i+1/\theta^2)}\left(\frac{t_1-1}{t_2-s}\right)^{i+1/\theta^2}\left(\frac{1-s}{t_2-s}\right)^k
.$$
\end{enumerate}

Thus $(Z_t)_{t\geq 0}$ is a well defined Markov process which determines Markov process $(X_t)_{t\geq 0}$ through the
 one-to-one transformation
 \rf{X2Z}.
\section{Conditional moments}\label{S: Moments}
We now verify that  $(X_t)_{t\geq 0}$  is a quadratic harness.
\begin{theorem} For $\theta>0$, let $(Z_t)$ be the Markov process defined in previous Section. Let $(X_t)$ be
defined by \rf{X2Z}. Then $(X_t)$ is the bi-Poisson process with
parameters $(\theta,\theta)$, i.e.  it has covariance \rf{cov},
conditional moments \rf{LH}, and \rf{q-Var} with $\eta=\theta$ and
$q=1$.
\end{theorem}
\begin{proof}
In Section \ref{S: Constr} we already verified that $\wo(X_t)=0$, $\wo{X_t^2}=t$.
Since $\calL(Z_t|Z_s)$ is binomial for $1\leq s<t$, we have $\wwo{Z_t}{Z_s}=\frac{s-1}{t-1}Z_s$.
Combining this with the already established martingale property for $t<1$,
we see   that $(X_t,\calF_{\leq t})_{t\geq 0}$ is a martingale.
From the martingale property we get \rf{cov}.

To compute the conditional moments, we calculate explicitly the conditional distribution of
$\calL(Z_t|Z_s,Z_u)$. These are routine calculations,
so we just state the final answers, and omit most of the calculations of the corresponding moments.
\begin{enumerate}
\item If $0<s<t<u<1$ then $\calL(Z_t-Z_s|Z_s,Z_u)$ is binomial with parameters $n=Z_u-Z_s$ and $p=\frac{(1-u)(t-s)}{(1-t)(u-s)}$.
Therefore
$$\wwo{Z_t}{\Fsu}=Z_s+\frac{(1-u)(t-s)}{(1-t)(u-s)}(Z_u-Z_s)$$
$$=\frac{(u-t)(1-s)}{(1-t)(u-s)}Z_s+\frac{(1-u)(t-s)}{(1-t)(u-s)}Z_u .$$
Using \rf{X2Z} we get
$$\wwo{X_t}{\Fsu}=-t/\theta+\frac{u-t}{u-s}(X_s+s/\theta)+\frac{t-s}{u-s}(X_u+u/\theta),$$ which gives
\rf{LH}.
Similarly,
$$\wVar{Z_t}{\Fsu}=\frac{(1-u)(t-s)(u-t)(1-s)}{(1-t)^2(u-s)^2}(Z_u-Z_s)$$ which
gives
$$\wVar{X_t}{\Fsu}=\frac{(t-s)(u-t)}{(u-s)^2}\big(\theta(1-s)(X_u+u/\theta)-\theta(1-u)(X_s+s/\theta)\big).$$
A calculation gives \rf{q-Var}.

\item If $0<s<t<1<u$ then $\calL(Z_t-Z_s|Z_s,Z_u)$ is negative binomial with parameters $r=Z_s+Z_u+1/\theta^2$ and
$p=\frac{(1-t)(u-s)}{(1-s)(u-t)}$. Therefore
$$\wwo{Z_t}{\Fsu}=
Z_s+r(1-p)/p$$
$$=\frac{(1-s)(u-t)}{(1-t)(u-s)}Z_s+\frac{(u-1)(t-s)}{(1-t)(u-s)}Z_u+\frac{(u-1)(t-s)}{\theta^2(1-t)(u-s)},$$ which leads to \rf{LH}
and
$$
\wVar{Z_t}{\Fsu}=\frac{r(1-p)}{p^2}$$
$$=\frac{(u-1)(1-s)(t-s)(u-t)}{(1-t)^2(u-s)^2}(Z_s+Z_u+1/\theta^2)
,$$ which after a calculation leads to \rf{q-Var}.

\item If $0<s<1<t<u$ then $\calL(Z_t-Z_u|Z_s,Z_u)$ is negative binomial with parameters $r=Z_s+Z_u+1/\theta^2$ and
$p=\frac{(t-1)(u-s)}{(t-s)(u-1)}$. A calculation verifies \rf{LH} and \rf{q-Var}.

\item If $1<s<t<u$ then $\calL(Z_t-Z_u|Z_s,Z_u)$ is binomial with $n=Z_s-Z_u$ and $p=\frac{(s-1)(u-t)}{(t-1)(u-s)}$.  A calculation verifies \rf{LH} and \rf{q-Var}.

\end{enumerate}
The conditional moments for the remaining choices of $s<t<u$ follow by continuity.
\end{proof}

\section{Additional Properties}\label{S: AddProp}

\begin{proposition}[Poisson representation]\label{two Poissons}
Let $(N_t^{\la})$ and $(M_t^{\la})$ be two independent
 Poisson processes  with intensity $\la>0$. If $(X_t)$ is a bi-Poisson process with positive
parameters $(\eta,\theta)$ then
 $$\calL\left.\left(\left(t\left(h(t)X_{\frac{\theta}{\eta
 h(t)}}+\frac{1}{\eta}\right)\right)_{t>
 0},\left(t\left(X_{\frac{\theta h(t)}{\eta }}+\frac{1}{\eta}\right)\right)_{t> 0}
 \right|X_{\theta/\eta}=\la-\frac{1}{\eta}\right)=\calL\left((N_t^{\la})_{t>
 0},(M_t^{\la})_{t>
 0}\right),
$$
where
$$h(t)=\frac{1+\theta t}{\theta t}.
$$
\end{proposition}
\begin{proof} Without loss of generality we assume $\eta=\theta$. By
\rf{X2Z} and the Markov property, it suffices to prove that
\begin{equation}\label{Z2P1}
\calL\left.\left(\left(Z_{1/h(t)}\right)_{t>0}\right|Z_1=\la\right)
= \calL\left((N_t^{\la})_{t>
 0}\right),
\end{equation}
and
\begin{equation}\label{Z2P2}
\calL\left.\left(\left(Z_{h(t)}\right)_{t>0}\right|Z_1=\la\right) =
\calL\left((M_t^{\la})_{t>
 0}\right).
\end{equation}
Both equalities follow now from elementary calculations of finite
dimensional distributions using the conditional distributions
identified in Section \ref{S: Constr}.

To prove \rf{Z2P2},  take $t_n<t_{n-1}<\dots t_1$ so that
$1<h(t_1)<h(t_2)<\dots h(t_n)$. Then for $k_1\geq k_2\geq \dots\geq
k_n$ denoting $Y_j=Z_{h(t_j)}$ we have
$$\Pr(Y_1=k_1, Y_2=k_2,\dots,Y_n=k_n|Z_1=\la)$$
$$
=\Pr(Y_n=k_n|Y_{n-1}=k_{n-1})\dots
\Pr(Y_2=k_2|Y_1=k_1)\Pr(Y_1=k_1|Z_1=\la)$$ $$=
\frac{\la^{k_1}\exp(-\la t_1)}{k_n!(k_{n-1}-k_n)!\dots
(k_1-k_2)!}t_n^{k_n}(t_{n-1}-t_n)^{k_{n-1}-k_n}\dots
(t_{1}-t_2)^{k_{1}-k_2},
$$
which proves \rf{Z2P2}.

The proof of \rf{Z2P1} is similar after using the generalized Bayes
formula: for $t_1<t_{2}<\dots t_n$ so that
$0<1/h(t_1)<1/h(t_2)<\dots 1/h(t_n)<1$. Then for $k_1\leq k_2\leq
\dots\leq k_n$,  denoting $Y_j=Z_{1/h(t_j)}$ we have

$$\Pr(Y_1=k_1, Y_2=k_2,\dots,Y_n=k_n|Z_1=\la)$$
$$=\frac{f_{Z_1|Y_n=k_n}(\la)}{f_{Z_1}(\la)}
\Pr(Y_n=k_n|Y_{n-1}=k_{n-1})\dots \Pr(Y_2 =k_2|Y_1
=k_1)\Pr(Y_1=k_1),
$$
where $f_{Z_1|Y_n=k_n}$ is the conditional density of $Z_1$ given
$Y_n=k_n$, and $f_{Z_1}$ is the density of $Z_1$, which are both
gamma, see the last paragraph of Section \ref{The pure birth phase}.
Elementary calculations now prove \rf{Z2P1}.
\end{proof}
\begin{proposition}[Time-inversion]\label{Time Inversion}
If $(X_t)$ is a bi-Poisson process with parameters $(\theta,\theta,q)$,
then $(tX_{1/t})_{t>0}$ has the same distribution as $(X_t)_{t>0}$.
(Compare \cite{Gallardo-Yor04}, \cite{Lawi-05}.)
\end{proposition}
\begin{proof} This follows from the fact that $(tX_{1/t})$ satisfies \rf{cov}, \rf{LH}, and \rf{q-Var}, and hence by
\cite[Theorem 4.5 and Proposition 4.13]{Bryc-Matysiak-Wesolowski-04} is determined uniquely. For $q=1$, the conclusion
can also be derived directly from \rf{X2Z}  and
the fact that Markov process $(Z_t)_{t>0}$ has the same transition probabilities as $(Z_{1/t})_{t>0}$.
\end{proof}
\begin{proposition}[Distribution of upward jumps]
For a bi-Poisson process $(X_t)$ with parameters $(\theta,\theta)$,
where $\theta>0$, define
$$
\Gamma_i=\sup\{s\in[0,1):\:Z_s=i\}\;,\;\;\;i=0,1,\ldots,
$$
i.e. $\Gamma_i$ is the time of the $(i+1)$-th jump of the process
$(X_t)$ from the line $y=\theta(1- t)i-t/\theta$, $0\le t<1$, for
$i=0,1,\ldots$, see Fig. \ref{Fig3}. Then the joint density of the
random vector $(\Gamma_0,\Gamma_1,\ldots,\Gamma_k)$ is
\begin{equation}\label{GGG}
f_{(\Gamma_0,\Gamma_1,\ldots,\Gamma_k)}(s_0,s_1,\ldots,s_k)
=\frac{\Gamma\left(\frac{1}{\theta^2}+k+1\right)(1-s_k)^{\frac{1}{\theta^2}+k-1}}
{\Gamma\left(\frac{1}{\theta^2}\right)(1-s_0)^2(1-s_1)^2\dots(1-s_{k-1})^2}
\end{equation}
for $0\le s_0<s_1<\ldots<s_k<1$ (and $0$ otherwise).
\end{proposition}
\begin{proof}  Let $(M_t)_{t\geq 0}$ be a homogeneous  pure birth process with birth rates
 $\la_n=n+1/\theta^2$, $n=0,1,\dots$. It is well known that the sojourn times $\tau_j$ of $(M_t)$ in state $j$
  are exponential with parameter $j+1/\theta^2$, $j=0,1,\dots$.
 For $0\leq t<1$ we have $Z_t=M_{-\ln(1-t)}$, so $-\log(1-\Gamma_k)=\sum_{j=0}^k\tau_j$. Therefore
 $\tau_k=\ln(1-\Gamma_{k-1})-\ln(1-\Gamma_k)$ (here, we set $\Gamma_{-1}=0$).
Since the Jacobian of the transformation $s_j\mapsto \ln(1-s_{j-1})-\ln(1-s_j),j=0,\dots,k$
is $J(s_0,\dots,s_k)=\prod_{j=0}^k(1-s_j)^{-1}$, and $\tau_0,\tau_1,\dots,\tau_k$ are independent, the joint density of
$(\Gamma_0,\Gamma_1,\ldots,\Gamma_k)$  is
$$
J(s_0,\dots,s_k)\prod_{j=0}^k\left((j+1/\theta^2)\exp\left(-(j+1/\theta^2)\ln\frac{1-s_{j-1}}{1-s_j}\right)\right)
$$
which simplifies to \rf{GGG}.

Alternatively, we can use the fact that
$$
R=P(s_0<\Gamma_0<s_1<\Gamma_1<s_2<\ldots<s_{k-1}<\Gamma_{k-1}<s_k<\Gamma_k)$$
$$=P(Z_{s_0}=0,Z_{s_1}=1,\ldots,Z_{s_k}=k).$$
The formula
$$
f_{(\Gamma_0,\Gamma_1,\ldots,\Gamma_k)}(s_0,s_1,\ldots,s_k)
=(-1)^k\frac{\partial^{k+1}\:R}{\partial s_k\partial
s_{k-1}\ldots\partial s_1\partial s_0}
$$
yields \rf{GGG} after a calculation.
\end{proof}
\begin{proposition}[Distribution of downward jumps]
Consider
$$
\Delta_i=\inf\{t>1:\:Z_t=i\}\;,\;\;\;i=0,1,\dots,
$$
i.e. $\Delta_i$ is the time of entrance of the process $(X_t)$ onto
the line $y=\theta( t-1)i-1/\theta$, $t>1$, for
$i=0,1,\dots$. Then the joint density of the random vector
$(\Delta_0,\Delta_1\ldots,\Delta_k)$ is
$$
f_{(\Delta_0,\Delta_1\ldots,\Delta_k)}(t_k,\ldots,t_1,t_0)
=\frac{\Gamma\left(\frac{1}{\theta^2}+k+1\right)(t_k-1)^{\frac{1}{\theta^2}+k-1}}
{\Gamma\left(\frac{1}{\theta^2}\right)
t_k^{\frac{1}{\theta^2}+k+1}( t_0-1)^2(t_1-1)^2\dots(t_{k-1}-1)^2}
$$
for $1<t_k<\ldots<t_1<t_0$  (and $0$ otherwise).
\end{proposition}
\begin{proof} From time-invertibility of the process, $(\Delta_0,\Delta_1\ldots,\Delta_k)$ has the same distribution as
$(1/\Gamma_0,1/\Gamma_1,\ldots,1/\Gamma_k)$. Thus
$$f_{(\Delta_0,\Delta_1\ldots,\Delta_k)}(t_k,\ldots,t_1,t_0)=\frac{1}{t_0^2t_1^2\dots t_k^2} f_{(\Gamma_0,\Gamma_1,\ldots,\Gamma_k)}(1/t_0,1/t_1,\ldots,1/t_k),$$
which simplifies to the expression above.
\end{proof}

%
%
%

\begin{proposition}[Time to reach lower boundary]
  The time a bi-Poisson process $(X_t)$ with parameters $\eta=\theta>0$, $q=1$ reaches the horizontal line $-1/\theta$
on which it stays forever is finite but has infinite expectation.
\end{proposition}
\begin{proof}
The distribution of  $\Delta_0=\inf\{t>1: X_t=-1/\theta\}$
is a special case of the distribution of jumps,
but it is just as easy to derive it independently. Since $Z_t$ is negative binomial, for $t>1$ we have
$\Pr(\Delta_0>t)=1-\Pr(Z_t=0)=1-(1-1/t)^{1/\theta^2}$. 

From the inequalities $(1-x)^p\leq 1-px$ when $0<p<1,x\geq 0$ and
$(1-x)^p\leq (1-x)$ when $p>1, 0\leq x\leq 1$
we get
$$\wo(\Delta_0)=\int_1^\infty\left(1-\left(1-\frac{1}{t}\right)^{1/\theta^2}\right)dt \geq
  \displaystyle\min\{1,1/\theta^2\}\int_1^\infty \frac{dt}{ t }=\infty.
$$
\end{proof}
\begin{proposition}[Poisson  limit]
For $\eta\theta> 0$ let $(X_t^{(\eta,\theta)})$ be the bi-Poisson
process with parameters $(\eta,\theta)$,  and let $(N_t)$ be the
Poisson process with parameter $\la=1$. As $\eta\to 0$ the process
$(\frac{1}{{\theta}}X_{t\theta^2}^{(\eta,\theta)})$ converges in
$D[0,\infty)$  to the
Poisson-type process $( N_t-t)_{t\geq 0}$.
\end{proposition}
\begin{proof}
Calculating the conditional variances one can check that
$(\frac{1}{{\theta}}X_{t\theta^2}^{(\eta,\theta)})$ is a bi-Poisson
process with parameters $(\eta\theta,1)$

Consider now the bi-Poisson process $(X_t)=(X_t^{(\eps)})$  with
parameters $(\eps,\eps)$ for $\eps=\sqrt{\eta\theta}$. Then by the
previous argument, process
 $Y^{\eps}=(\frac{1}{{\eps}}X_{t\eps^2}^{(\eps)})_{t\geq 0}$ has the
 same distribution as process
 $(\frac{1}{{\theta}}X_{t\theta^2}^{(\eta,\theta)})$.
Therefore it suffices to show that as $\eps\to0$,  the process
$Y^{\eps}$ converges in $D[0,\infty)$  to the
Poisson-type process
$( N_t-t)_{t\geq 0}$.

 We first verify the convergence of finite-dimensional distributions.
For $0\leq t<1/\eps^2$, the appropriate version of \rf{X2Z} is
$$\frac{1}{{\eps}}X_{t\eps^2}^{(\eps)}=(1-t\eps^2)Z_{t\eps^2}^{(\eps)}-t,$$
where $Z_{t\eps^2}^{(\eps)}$ is a (non-homogeneous) pure birth
process  on $0\leq t<1/\eps^2$. We will verify that the finite
dimensional distributions of $(Z_{t\eps^2}^{(\eps)})$ converge to
the finite dimensional distributions of $(N_t)$.

Fix arbitrary $0=t_0<t_1<t_2<\dots<t_n$ and
$u_1,u_2,\dots,u_n\leq 0$. It suffices to show that

\begin{equation}\label{PL}
\lim_{\eps\to 0}\wo\left(\exp\left(\sum_{j=1}^n u_j
Z_{t_j\eps^2}^{(\eps)} \right)\right)=\prod_{j=1}^n
\exp\left(-(t_j-t_{j-1})(e^{\sum_{i=j}^nu_i}-1)\right)
\end{equation}

We rely on the following observation, which can be regarded as
special case of Slutsky's theorem: if $u_\alpha\to u$ and
$(W_{1}^{(\alpha)}, W_{2}^{(\alpha)},\dots,W_{n}^{(\alpha)})$
converges weakly to a random vector $(W_{1},W_{2},\dots,W_{n})$ as
$\alpha\to 0$ and appropriate exponential moments exist, then
\begin{equation}\label{Slucki}
\lim_{\alpha\to 0}\wo\left(\exp\left(\sum_{j=1}^n u_j
W_{j}^{(\alpha)}+u_\alpha W_{n}^{(\alpha)}\right)\right) =
\wo\left(\exp\left(\sum_{j=1}^n u_j W_{j}+u W_{n}\right)\right).
\end{equation}

 From \rf{MGF} we have
$$\wo\left(\exp\left(\sum_{j=1}^n u_j Z_{t_j\eps^2}^{(\eps)}
\right)\right)=A_\eps \wo\left(\exp\left(\sum_{j=1}^{n-1} u_j
Z_{t_j\eps^2}^{(\eps)} +(u_n+u_\eps)
Z_{t_{n-1}\eps^2}^{(\eps)}\right)\right),$$ where
$$e^{u_\eps}=\frac{1-\eps^2
t_n}{1-\eps^2 t_{n-1}-\eps^2(t_n-t_{n-1})e^{u_n}}\to 1
$$
and $$ A_\eps=e^{u_\eps/\eps^2}\to e^{-(t_n-t_{n-1})(e^{u_n}-1)}.$$
This proves \rf{PL} for $n=1$, and shows that \rf{PL} holding for
$n-1$ implies \rf{PL} for $n$, ending the proof by induction.
%
%

 Therefore, the increments of $(Z_{t\eps^2}^{(\eps)})$ are
asymptotically independent Poisson-distributed with parameter
$\la=t-s$, and
 the  finite-dimensional distributions of $(Z_{t\eps^2}^{(\eps)})$ converge to the corresponding distributions of $(N_t)$.

Tightness of $Y^\eps$ in $D[0,\infty)$ follows from \cite[Theorem 2.1 and Remark 4.2]{Gut-Janson-01}.
Indeed, their condition (i$^\infty$) holds as $Y^\eps$ is a martingale and $E(Y^\eps(t))^2=t$. Their
condition (ii$^\infty$) holds as for fixed $a_1\leq s \leq a_2$ and all $\eps$ small enough from \rf{MGF} we get
$$
\Pr(\mbox{ at least two $Y^\eps$ jumps in } [s,s+\delta)) $$ $$=
\Pr(X^{(\eps)} \mbox{ has at least two jumps in }
[s\eps^2,(s+\delta)\eps^2))$$
$$
\leq \Pr( Z_{(s+\delta)\eps^2}^{(\eps)}-Z_{s\eps^2}^{(\eps)}\ge
2)\leq \delta^2 C(a_1,a_2) E(1+\eps^2Z_{a_2}^{(\eps)})^2.
$$
Thus $$\delta^{-1}\limsup_{\eps\to 0}\sup_{a_1\leq s\leq a_2}\Pr(\mbox{ at least two $Y^\eps$ jumps in } [s,s+\delta))
\leq \delta C(a_1,a_2)\to 0 \mbox{ as } \delta\to 0.
$$
Finally, their condition (iii) holds as for fixed $\eta>0$ and all $a>0$ small enough
$$\limsup_{\eps\to 0}\Pr(\sup_{s\leq a} |Y^{\eps}(s)|>\eta)\leq \limsup_{\eps\to 0}\Pr(Z_a^{(\eps)}\geq 1)=ae^{-a}\to 0 \mbox{ as } a\to 0. $$
%
%
%
\end{proof}
\begin{proposition}[Brownian limit]
Let $(X_t)_{t\geq 0}=(X_t^{(\theta)})_{t\geq 0}$ be the bi-Poisson process with parameters $\eta=\theta$, $q=1$.
Then as $\theta\to0$ the process
$(X_t^{(\theta)})_{t\geq 0}$ converges in $D[0,\infty)$  to the
standard Brownian motion
$(B_t)_{t\geq 0}$.
\end{proposition}
\begin{proof}
The convergence of finite dimensional distribution follows from uniqueness of the quadratic harnesses (see
 \cite{Wesolowski93} for the special case we need here) 
 and the fact that
the limiting process must satisfy \rf{q-Var} by the uniform integrability of
$\{(X_t^{(\theta)})^2:\theta \leq 1\}$; in fact, $\wo(X_t^4)=2 t^2+ ( t + 3t^2 + t ^3 )\theta^2$.

Tightness in $D[0,\infty)$ now follows from \cite[Proposition 1.2]{Aldous-89}, as $(X_t^{(\theta)}:0\leq t<\infty)$
is a martingale for each $\theta$, and the limiting process is continuous.

Here we give  a simple direct argument for the convergence of finite
dimensional distributions. It is enough to prove that for all
$0=t_0<t_1<t_2<\dots<t_n$ and all
 $u_1,u_2,\dots,u_n$ close enough to zero
\begin{equation}\label{ME}
\lim_{\theta\to 0}\wo\left(\exp\left(\sum_{j=1}^n u_j
X_{t_j}^{(\theta)}\right)\right)= \exp\frac12\left(\sum_{k=1}^n
(t_k-t_{k-1})\left(\sum_{j=k}^n u_j\right)^2\right)
\end{equation}
We proceed by induction, suppressing $\theta$ in $X_{t}^{(\theta)}$
to shorten the expressions. Without loss
of generality we may assume $t_1<1$ (set $u_1=0$, if necessary). To verify \rf{ME} for $n=1$ and $t_1<1$ we
use \rf{X2Z} and \rf{MGF} with $s=0$, which gives
$$
\wo\left(e^{u_1X_{t_1}}\right)=\left(\frac{1-t_1}{e^{u_1 t_1\theta}
- t_1e^{u_1\theta}}\right)^{1/\theta^2}\to \exp\left({\frac12u_1^2
t_1}\right).
$$

Suppose \rf{ME} holds for $0=t_0<t_1<t_2<\dots<t_{n-1}<1$ and let
$t_n\in(t_{n-1},1)$. Then again using \rf{X2Z} and \rf{MGF} we get
\begin{equation}\label{AF}
\wo\left(\exp\left(\sum_{j=1}^n u_j X_{t_j}\right)\right)=
A_\theta\wo\left(\exp\left(\sum_{j=1}^{n-1} u_j X_{t_j}+u_\theta
X_{t_{n-1}}\right)\right),
\end{equation}
where
$$
A_\theta=\left(\frac{1-t_n}{1-t_{n-1}-(t_n-t_{n-1})e^{\theta
u_n(1-t_n)}}\right)^{{1}/{(\theta^2(1-t_{n-1}))}}\exp\left(\frac{-u_n(t_n-t_{n-1})}{\theta(1-t_{n-1})}\right),
$$
$$
e^{u_\theta}=\left(\frac{1-t_n}{1-t_{n-1}-(t_n-t_{n-1})e^{\theta
u_n(1-t_n)}}\right)^{{1}/{(\theta(1-t_{n-1}))}}\exp\left(u_n\frac{1-t_n}{1-t_{n-1}}\right)
$$
Since
$$
A_\theta=\left(\frac{1-t_n}{(1-t_{n-1})e^{\theta
u_n(t_n-t_{n-1})}-(t_n-t_{n-1})e^{\theta
u_n(1-t_{n-1})}}\right)^{{1}/{(\theta^2(1-t_{n-1}))}},
$$
by Taylor expansion $\lim_{\theta\to
0}A_\theta=e^{\frac12u_n^2(t_n-t_{n-1})}$. Similarly
$\lim_{\theta\to 0}u_\theta=u_n$, which using \rf{Slucki} and
induction assumption proves  \rf{ME}.

Now we know that \rf{ME} holds for all $n\geq 1$ and
$0=t_0<t_1<t_2<\dots<t_{n-1}<1$,  and we take $t_n=1$.
Using \rf{X2Z} and \rf{MG Gamma} we see that \rf{AF} holds with
$$
e^{u_\theta}=\left(\frac{1}{1-(1-t_{n-1})u_n\theta}\right)^{{1}/{(\theta(1-t_{n-1}))}}\to
e^{u_n}
$$
and
$$
A_\theta=e^{(u_\theta-u_n)/\theta}\to e^{\frac12
u_n^2(1-t_{n-1})}.
$$
By \rf{Slucki} this proves \rf{ME}.

Now we know that \rf{ME} holds for all
$0=t_0<t_1<t_2<\dots<t_{n-2}<t_{n-1}=1$, $n\geq 1$ and we take
$t_n>1$. By \rf{X2Z} and \rf{MG Poiss} we get \rf{AF} with
$$
e^{u_\theta}=\exp\left(\frac{e^{\theta u_n
(t_n-1)}-1}{\theta(t_n-1)}\right)\to e^{u_n}
$$
and $$A_\theta=e^{(u_\theta-u_n)/\theta}\to e^{\frac12
u_n^2(t_n-1)}.$$ Again using
 \rf{Slucki} we get \rf{ME}.

Finally, we assume that \rf{ME} holds for all
$0=t_0<t_1<t_2<\dots<t_{n-1}$, $n\geq 1$ with $t_{n-1}>1$ and we
take $t_n>t_{n-1}$.

Using \rf{X2Z} and \rf{MG Gamma} we see that \rf{AF} holds with
$$
e^{u_\theta}=\left(\frac{t_n-t_{n-1}+(t_{n-1}-1)e^{\theta u_n
t_n-1}}{t_n-1}\right)^{{1}/{(\theta(t_{n-1}-1)})}\to e^{u_n}
$$
and
$$
A_\theta=e^{(u_\theta-u_n)/\theta}\to e^{\frac12
u_n^2(t_n-t_{n-1})}.
$$
By \rf{Slucki} this proves \rf{ME}.

\end{proof}

\subsection*{Acknowledgement} We would like to thank
Wojciech Matysiak for comments and discussions and Sergiusz
Weso{\l}owski for assistance with
simulations of the bi-Poisson process in Fig. \ref{Fig3}. 

\bibliographystyle{plain}
\bibliography{bipois,Vita,Wesol,q-reg,q-nc,harness,free-lalu}
\end{document}